\documentclass[12pt]{article}

\usepackage{epsfig}

\newtheorem{theorem}{Theorem}[section]
 \newtheorem{lemma}[theorem]{Lemma}
\newtheorem{definition}[theorem]{Definition}

 \newtheorem{proposition}[theorem]{Proposition}
\newenvironment{proof}{{\it Proof:\/}}{$\Box$\vskip 0.08in}
\usepackage{amssymb}
\begin{document}
\begin{center} 
{\large\bf The topological interpretation of
the core group of a
surface in $S^4$} 
\end{center}
 \begin{center}
by J\'ozef H.Przytycki\footnote{Supported by NSF grant:
NSF: DMS--9808955}
 and Witold Rosicki\footnote{Supported by UG grant: BW-5100-5-0256-9}
 \end{center} 
 \vspace{0.1in}
\begin{quotation}\ \\  
{\bf Abstract.}
We give a topological interpretation of the core group 
invariant of a surface embedded in $S^4$ \cite{F-R,Ro}.
We show that the group is isomorphic to the free product of
the fundamental group of the double branch cover of $S^4$
with the surface as a branched set, and the infinite cyclic group.
We present a generalization for unoriented surfaces, for other
cyclic branched covers, and other codimension two embeddings of
 manifolds in spheres.
\end{quotation}

It is shown in \cite{Ya-1} (compare \cite{Ka-1,C-S-2,Ya-2}) that the 
fundamental group of the complement of an oriented surface, $M$,
in $S^4$ allows a Wirtinger type presentation.
In \cite{Ro} the core group invariant of $M$ was introduced
(following \cite{F-R,Joy}) and the topological interpretation
of the group was promised. In the first section we give
 this interpretation
following the result of Wada's \cite{Wa} for classical knots
and the proof of Wada result presented in \cite{Pr}.   
In the second section we generalize our results to
unoriented surfaces, and codimension two embeddings of
closed $n$-manifolds in $S^{n+2}$.

\section{Surfaces in $S^4$}\label{1}

Let $M$ be an oriented surface embedded in $S^4= R^4 \cup\{ \infty\}$.
Let $p:R^4 \to R^3$ be a projection such that the restriction to $M$,
$p|_M$ is a general position map. 
Define the lower decker set (``invisible" set) \cite{C-S-1,C-S-2}
as $A=\{ x\in M  \ | \ \exists y\in M,\ p(x)=p(y)\ and \ \ x \
is\ the\ lower\ point,\ that\ is\ x \  is\ further \ than \ y \\ in \ the \ 
direction \ \ of \ \ the \ \ projection \}$. The set $A$ separates $M$ into
regions (``visibility" regions), that is arc connected components 
of $M - A$.
\begin{theorem}[\cite{Ya-1}]
The fundamental group of $S^4-M$ ($\Pi_M =\pi_1(S^4 - M)$) 
has the following Wirtinger type presentation. 
Generators of the group, $x_1,x_2,...,x_n$,
correspond to regions of $M$. Relations of the group correspond 
to double point arcs and are of the form $x_ix_px_i^{-1}x_q^{-1}$ or
$x_ix_p^{-1}x_i^{-1}x_q$, depending on the orientations of $p(M)$ at
double point arcs.
Here $x_i$ corresponds to the higher region of the projection and $x_p,
x_q$ to lower regions.
\end{theorem}
Similarly to the group $\Pi_M$, 
we define the core group, $G_M$ of a surface $M$ in $S^4$.

\begin{definition}[\cite{Ro}]\label{1.2}
The core group, $G_M$ of a surface $M$ in $S^4$, has the
following, quandle type, presentation. 
Generators of the group, $y_1,y_2,...,y_n$, correspond to regions of $M$.
Relations of the group correspond
to double point arcs and are of the form $y_iy_p^{-1}y_iy_q^{-1}$.
Here $y_i$ corresponds to the higher region of the projection and $y_p,
y_q$ to lower regions.
\end{definition}

The operation $y_p*y_i = y_iy_p^{-1}y_i$ defines a quandle.\\
We will interpret topologically the core group as follows.
It is the analogue of the Wada result for classical links \cite{Wa}.

\begin{theorem}\label{1.3}
The core group, $G_M$ of a surface $M$ in $S^4$, is isomorphic to the
free product of the fundamental group of the branched double cover of
$S^4$ with branching set $M$, and the infinite cyclic group.
That is $G_M= \pi_1(M^{(2)}) *Z$ where $M^{(2)}$ is the considered
double branched cover.
\end{theorem}

Below we will prove a generalization 
of the theorem for any cyclic 
branched cover, 
following the exposition of the
classical case in \cite{Pr}.

Let $f: M^{(k)} \to S^4$ be a cyclic $k$ fold branched covering 
with a branch set $M$. More precisely $M^{(k)}$ is defined as follows:
Let $V_M$ be a tubular neighborhood of $M$ in $S^4$. It has a 
structure of a 2-disc bundle (over $M$). For a point $x \in M$
the boundary of the disk $D_x$ of the bundle is called a meridian
of $M$ (at $x$), denoted by $\mu_x$. 
It has the natural orientation for oriented $M$.
The first homology group of $S^4 - M$ is freely generated by
meridians of $M$ (one for each connected component of $M$).
Thus, we have a k-fold cyclic covering of $S^4 - M$  given
by an epimorphism $H_1(S^4 - M) \to Z_k$ which sends $\mu_x$ to 1.
This covering can be uniquely extended to the cyclic $k$ 
fold branched covering
with a branch set $M$. Namely, we define the cover on each (meridian) disk
$D_x$ as given by the function $p: D'_x \to D_x$,
$p(z) = z^k$ where $z$ is a complex coordinate of the disks.

Let $\Pi_M^k = \pi_1(M^{(k)})$.
\begin{proposition}
Consider the epimorphism $\hat g_k: \pi_1(S^4 - M) \to Z_k$
given by $\hat g_k (x_i) = 1$ for any $i$. Then
$$ \Pi_M^k = ker (\hat g_k)/(x_i)^k.$$
\end{proposition}
\begin{proof}
Consider a regular neighborhood of $M$ in $S^4$ and 
its normal disk bundle. 
Each disk in the bundle is covered by a disk 
( the center being a branch point). As noted before the covering
 can be described 
by the map $z\to z^k$ where $z$ is a complex coordinate of the disk. 
In the Wirtinger presentation of the fundamental group,
generators correspond to some meridians (that is, boundaries of
disks of the normal bundle).\\
Consider the unbranched k-covering of $S^4-M$ obtained by removing
the branch set.
The epimorphism $\hat g:\pi_1(S^4-M)\to Z_k$ corresponding to 
the covering sends each meridian to $1 \in Z_k$, ($\hat g(\mu_x)= 1$).
The fundamental group of the unbranched cover is equal to 
$ker \hat g$.
To get a branched cover from the unbranched one we have to 
fill the cover of every meridian. On algebraic level we add 
relations $\mu_x^k = 0$.
Because every meridian is conjugated to a generator $x_i$,
the proposition is proved.
\end{proof}

\begin{theorem}\label{1.5} 
$\Pi_M^{(k)}*\underbrace{Z*...*Z}_{k-1}= \Pi_{M\sqcup S^2}^{(k)}$
has the following quandle type description (in the disjoint sum $M\sqcup S^2$,
 $S^2$ is the unknotted component that can be separated by a $3$-sphere 
from $M$).\\
To every (visible) region there correspond $k-1$ generators 
$\tau^j(y_i)$  ($0 \leq j<k-1$) and to each double point arc correspond 
$k-1$ relations $\tau^j(y_p\tau(y_i)\tau(y_q^{-1})y_i^{-1})$  or 
$\tau^j(y_i\tau(y_p)\tau(y_i^{-1})y_q^{-1})$  depending on 
a local orientation.
\end{theorem}

\begin{lemma}\label{1.6}
 Let $F_{n+1} = \{x_1,x_2,...,x_{n+1}\  | \ \}$ 
be a free group of $n+1$  generators, $x_1,x_2,...,x_{n+1}$. 
Consider the epimorphism $g_{k}: F_{n+1} \to Z_k$ 
such that $g_k(x_i)= 1$ for any $i$. Let $F^{(k)}= ker(g_k)$,
 and let $\bar F^{(k)}=F^{(k)}/(x_i^k)$.
Let $\tau :F_{n+1}\to F_{n+1}$ be an automorphism, defined by
 $\tau(x_i)=x_{n+1} x_i x_{n+1}^{-1} $ and $y_i=x_i x_{n+1}^{-1}$.
Then
\begin{enumerate}
\item[(i)]
 $F^{(k)}$ is a free group of $nk+1$  generators : $x_{n+1}^k$ and 
$\tau^j(y_i)$ for $0 \leq j\leq k-1$, $1\leq i\leq n$.\\
\item[(ii)]
 $\bar F^{(k)}$  is a free group on 
$n(k-1)$ generators $\tau^j(y_i)$ 
for $0\leq j < k-1$ , $1\leq i\leq n$.\\
Furthermore $y_i \tau(y_i)...\tau^{k-1}(y_i) = 1$ for any $i$.\\
\end{enumerate}
\end{lemma}
\centerline{\psfig{figure=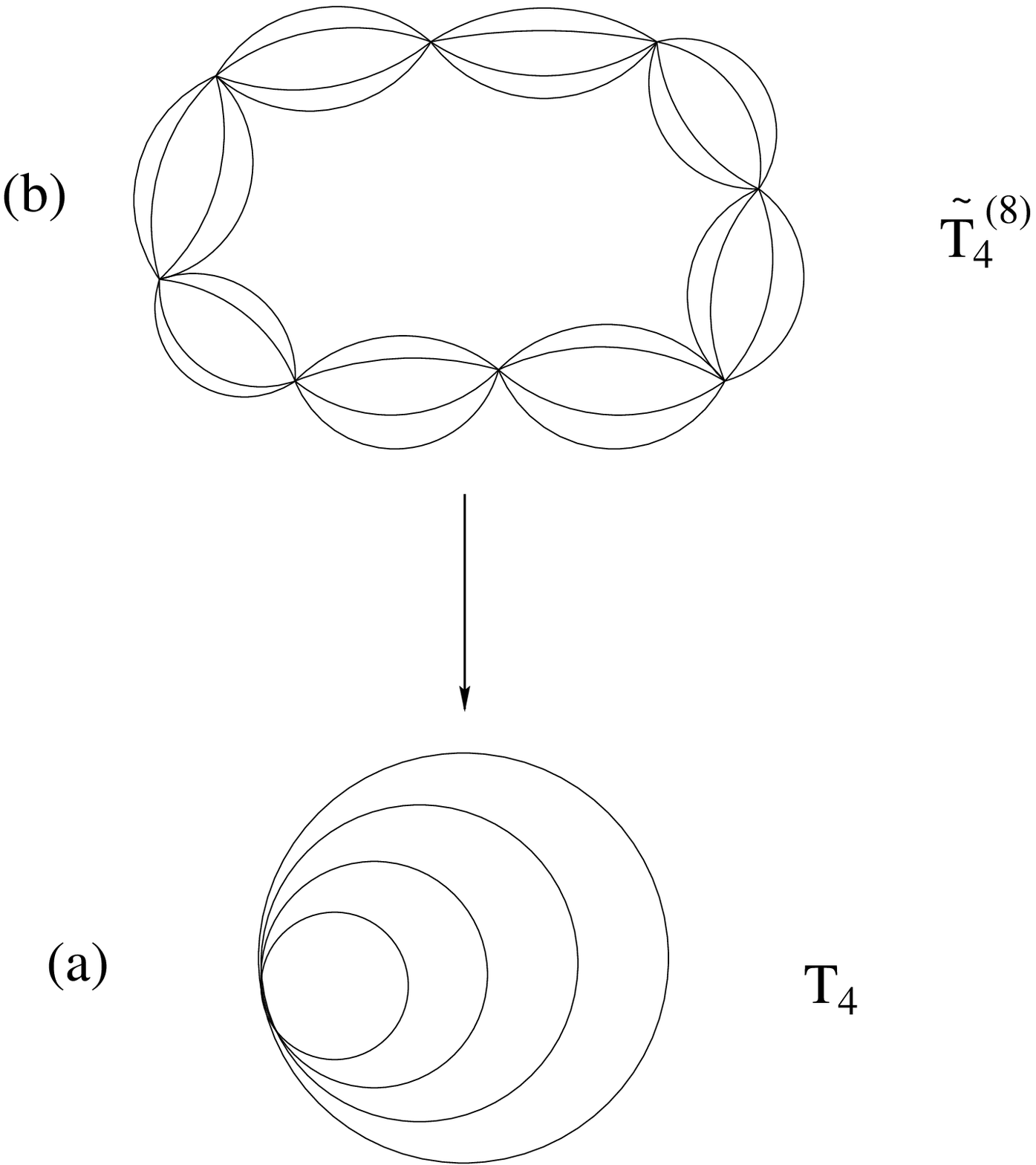,height=7.2cm}}
\begin{center}
Fig. 1  \ \ (The covering $\tilde T_4^{(8)} \to T_4$)
\end{center}
\begin{proof}
\begin{enumerate}
\item[(i)] 
$F_{n+1}= \pi_1(T_{n+1})$ where $T_{n+1}$ is a graph of one vertex and 
$n+1$ loops (Fig.1(a)).  $g_k:F_{n+1}\to Z_k$ determines a k-fold
 covering space
$\tilde T_{n+1}^{(k)}$ of the graph $T_{n+1}$ as show in Fig.1(b).
Of course, the Euler characteristic $\chi (T_{n+1})=-n$ and
 $\chi(\tilde T_{n+1}^{(k)})=-kn$. Thus, $\pi_1(\tilde T_{n+1}^{(k)})$
 is freely
generated by $nk+1$ generators. We can easily identify free generators
 of $\pi_1(\tilde T_{n+1}^{(k)})= ker g_{k}$; in particular, 
$g_{k}(x_{n+1}^k) = g_{k}(\tau^j(y_i)) = 0$ 
for $1\leq i\leq n$ and $0\leq j\leq k-1$.
\item[(ii)] 
  $y_i\tau (y_i)\tau^2(y_i)...\tau^{k-1}(y_i)= 
(x_ix_{n+1}^{-1})(x_{n+1}x_ix_{n+1}^{-2})...(x_{n+1}^{k-1}x_ix_{n+1}^{-k})= 
x_i^kx_{n+1}^{-k}$.\\
Thus, in $F^{(k)}_{n+1}/(x_{n+1}^k)$ we have $y_1...\tau^{k-1}(y_i)= x_i^k$
and furthermore\\
$\bar F_{n+1}^{(k)} = 
\{F_{n+1}^{(k)}\ | \ x_{n+1}^k, x_i^k\ ; \ 
 1\leq i \leq n \}$=\\
$\{x_{n+1}^k, \tau^j(y_i) \ | \  x_{n+1}^k, y_i \tau(y_i)...\tau^{k-1}(y_i) 
,1\leq i \leq n , 0\leq j\leq k-1\} =$ \\
$=\{ \tau^j(y_i) \ 
| \  y_i \tau (y_i)...\tau^{k-1}(y_i),
1\leq i \leq n , 0\leq j\leq k-1 \} =
\{ \tau^j(y_i) \ | \  1\leq i \leq n , 0\leq j <k-1 \}$  
as required.
\end{enumerate}
\end{proof}

We will use an extension of Lemma 1.6 for a group with Wirtinger 
type relations.\\
Let $G=\{F_n|\ r_1,...,r_m\}=\{x_1,...,x_n|\ r_1,...,r_m \}$ 
where any relation
 is of the form $r=x_ix_p^{\epsilon}x_i^{-1}x_q^{-\epsilon}$ 
for $\epsilon =\pm 1$, and $1\leq i,p,q \leq n$\\
We have $G*Z = \{x_1,...,x_{n+1}| r_1,...,r_m\}$. The epimorphism 
$g_k: F_{n+1} \to Z_k$ yields the epimorphism
$ \hat g_k: G*Z \to  Z_k$ . The epimorphism 
$\hat g_k$ is well defined because $g_k(r_i)=0$.
We use Lemma 1.6 to find a presentation of  $G^{(k)} = ker \hat g_k$ 
and $\bar G^{(k)} = G^{(k)}/(x_i^k)$.\\
\begin{lemma}\label{1.7}
\begin{enumerate}
\item[(i)] $G^{(k)} = 
\{F^{(k)}| \tau^j (r_s), 1\leq s \leq m , 0\leq j \leq k-1 \}.$
\item[(ii)] 
$\bar G^{(k)} = \{\bar F^{(k)} \ | \  
\tau^j (r_s), 1\leq s \leq m , 0\leq j < k-1 \}
= \{\tau^j(y_i) | \ \tau^j (r_s)\ , 1\leq s \leq m , 0\leq j < k-1 ,
 1\leq i \leq n \}$  where \ 
$\tau^j (x_ix_px_i^{-1}x_q^{-1}) = \\ 
\tau^j (y_i \tau (y_p) \tau (y_i^{-1})y_q^{-1})$
and
$\tau^j (x_ix_p^{-1}x_i^{-1}x_q) = 
\tau^j(y_iy_p^{-1}\tau^{-1}(y_i^{-1})\tau^{-1}(y_q))$.
\end{enumerate}
\end{lemma}
\begin{proof}
\begin{enumerate}
\item[(i)]
 Because $r_s \in ker g_k$ then the relations of 
$ker \hat g_k$ are of the form 
$wr_iw^{-1}$ for $w \in F_{n+1}$. We observe that
if $g_k(w)=u$ ($0\leq u \leq k-1$) then  $w=x_{n+1}^uw'$
where $w'\in ker g_k$.
Therefore relations of $G^{(k)}= ker \hat g_k$ are of the form 
$\tau^u (r_s)$ for $0\leq u \leq k-1 , 1\leq s \leq m$. 
These yield the presentation of $G^{(k)}$. We have:\ \
$x_ix_px_i^{-1}x_q^{-1}= \\
(x_ix_{n+1}^{-1})(x_{n+1}(x_px_{n+1}^{-1})
x_{n+1}^{-1})(x_{n+1}(x_{n+1}x_i^{-1})x_{n+1}^{-1})(x_{n+1}x_q^{-1})=\\  
 y_i \tau (y_p) \tau (y_i^{-1})y_q^{-1}$\\
$x_ix_p^{-1}x_i^{-1}x_q =\\ 
(x_ix_{n+1}^{-1})(x_{n+1}x_p^{-1})(x_i^{-1}x_{n+1})(x_{n+1}^{-1}x_q)=
y_iy_p^{-1}(\tau^{-1}(y_i))^{-1} \tau^{-1}(y_q)$

\item[(ii)] Adding relations $(x_i)^k$ reduces $F^{(k)}$ to 
$\bar F^{(k)}$ and $G^{(k)}$ to $\bar G^{(k)}$ so
$\bar G^{(k)}=
\{\bar F^{(k)}|\ \tau^j(r_s), 1\leq s\leq m, 0\leq j\leq k-1 \}$
We can eliminate the relation $\tau^{k-1}(r_s)$ 
(expressing it using other relations).
We use the identity\\
$y_i\tau (y_i)...\tau^{k-1}(y_i)=1$ in $\bar G^{(k)}$
(in particular $\tau^{k}(y_i)=y_i$).  We can write our
relations $\tau^u(y_i\tau(y_p)\tau(y_i^{-1})y_q^{-1})$ as 
$\tau^{u+1}(y_py_i^{-1})
= \tau^u(y_i^{-1}y_q)$. We assume it holds for $0 \leq u <k-1$ and we will
see that $\tau^{k-1}(r_s)=1$ as well:
$$\tau^{k-1}(y_i\tau(y_p)\tau(y_i^{-1})y_q^{-1})= \tau^{k-1}(y_i)\tau^k(y_p)
\tau^k(y_i^{-1})\tau^{k-1}(y_q^{-1})= $$
$$=\tau^{k-1}(y_i)y_py_i^{-1}\tau^{k-1}(y_q^{-1}) =
\tau^{k-1}(y_i)y_p(y_i^{-1}y_q)\tau(y_q)...\tau^{k-2}(y_q) =$$
 $$= \tau^{k-1}(y_i)y_p\tau(y_p)(\tau(y_i^{-1})\tau(y_q))\tau^2(y_q)
...\tau^{k-2}(y_q) =$$  
$$=\tau^{k-1}(y_i)y_p\tau(y_p)\tau^2(y_p)
(\tau^2(y_i^{-1})\tau^2(y_q))\tau^3(y_q)...\tau^{k-2}(y_q)= $$
$$= \tau^{k-1}(y_i)y_p\tau(y_p)\tau^2(y_p)\tau^3(y_p)...\tau^{k-1}(y_p)
\tau^{k-1}(y_i^{-1}) = 1$$ as required.
\end{enumerate}
\end{proof}
Theorem 1.5 follows from Proposition 1.4 and Lemma 1.7.
In particular, Theorem 1.3 follows because for $k=2$ we have $y_i \tau (y_i)= 
1$ so $\tau (y_i) = y_i^{-1}$. 
In the next section, we show that the recent result of Kamada \cite{Ka-2}
allows us to extend Theorem 1.5 to any closed oriented $n$-manifold 
in $S^{n+2}$. Similarly Theorem 1.3 can be generalized to any closed
unoriented $n$-manifold in $S^{n+2}$.

\section{Higher dimensional case: $M^{n} \subset S^{n+2}$.}\label{2}
We can extend our results to the case of
a closed n-manifold $M$ in $S^{n+2}$ (not necessary oriented or connected). 
First we need an existence of the Wirtinger type presentation
of the fundamental group of $S^{n+2}-M$  for oriented M, 
described in \cite{Ka-2} and its variant for unoriented $M$. 
\begin{theorem}\label{2.1}([Kamada])\ \\
Let $p: R^{n+2} \to R^{n+1}$ be a projection and $M$ be in general
position with respect to the projection. We can think that the base point
is very high above $R^{n+1}$ (say at $\infty \in S^{n+2}$). 
We have $n-1$ dimensional strata being 
the closure of the (invisible) set $\{x\in M|\
\exists y\in M,\ p(x)=p(y)\ and \ \ x \
is\ the\ lower\ point\}$. These strata cut $M$ into n-dimensional 
regions (``visibility regions").
Let $\Pi_M $ denote the fundamental group of $S^{n+2}-M$ 
($\Pi_M =\pi_1(S^{n+2} - M)$).
\begin{enumerate}
\item[(1)] Assume that $M$ is oriented.\\
Then $\Pi_M$ has the following Wirtinger type presentation.
Generators of the group, $x_1,x_2,...,x_n$,
correspond to regions of $M$ (and can be visualize by joining
the base point $b$ to any meridian of the region; meridians
are oriented so elements of $\Pi_M$ are well defined).
Relations of the group correspond
to double point strata and are of the 
form $x_ix_px_i^{-1}x_q^{-1}$ or $x_ix_p^{-1}x_i^{-1}x_q$, 
depending on orientations of $p(M)$ at double point strata.
Here $x_i$ corresponds to the higher region of the projection and $x_p,
x_q$ to lower regions.
\item[(2)] Consider now $M$ not necessary oriented. Let
generators, $x_1,x_2,...,x_n$ be chosen as before (one for each
 region). Because a region may be unorientable 
we have to make choices here. However, modulo relations $x_i^2=1$ 
(for $1\leq i\leq n$), the group has Wirtinger type presentation.
That is: $\Pi_M/(x_i^2)$ has a presentation as in (1) with
additional relations $x_i^2=1$.
\end{enumerate}
\end{theorem}
In the unoriented case meridians have no preferred orientation
so only the branched double cover is uniquely defined, so we can
interpret the core group of the embedded unoriented codimension
two submanifold of $S^{n+2}$.

Kamada's theorem allows us to generalize Theorems 1.3 and 1.5.
The proof, as before, bases on Lemma 1.7 ($k=2$ yields an unoriented case),
so we omit it.

Let $p: R^{n+2} \to R^{n+1}$ be a projection and $M$ be an
$n$-dimensional closed submanifold of $R^{n+2}$ 
(and $S^{n+2} = R^{n+2} \cup \infty$) being in general
position with respect to the projection. Using notation of
Theorem 2.1 we have.
\begin{theorem}\label{2.2}
Let $M^{(k)}$ denote the cyclic branched regular $k$-covering of
$S^{n+2}$ with an oriented manifold $M$ as a branched set, and
$\Pi_M^{(k)}$ its fundamental group. Then
$\Pi_M^{(k)}*\underbrace{Z*...*Z}_{k-1}= \Pi_{M\sqcup S^n}^{(k)}$
has the following quandle type description.\\
To every region corresponds $k-1$ generators
$\tau^j(y_i)$  ($0 \leq j<k-1$) and to each double point strata correspond
$k-1$ relations $\tau^j(y_p\tau(y_i)\tau(y_q^{-1})y_i^{-1})$  or
$\tau^j(y_i\tau(y_p)\tau(y_i^{-1})y_q^{-1})$  depending on a local orientation.
\end{theorem}
\begin{theorem}\label{2.3}
Let $M^{(2)}$ denote the double branched regular $2$-covering of
$S^{n+2}$ with an unoriented (possibly not orientable)
 manifold $M$ as a branched set. Then
$\pi_1(M^{(2)}) *Z = \Pi_{M\sqcup S^n}^{(2)}$ 
has the following core (quandle type) description.\\ 
Generators of the group, $y_1,y_2,...,y_n$, correspond to 
regions of $M$.
Relations of the group correspond
to double point strata and are of the form $y_iy_p^{-1}y_iy_q^{-1}$.
\end{theorem}
In the case of a 3-manifold in $S^5$, one can check invariantness
 of our groups,
defined combinatorically, using the Roseman moves \cite{Ros-1}.
\section{Acknowledgment.} 
We would like to thank Prof. L. Kauffman
for his encouragement to work on the above problems after
the talk of the second author at the {\it Knots in Hellas 98} conference.

\ \\
J\'ozef H.Przytycki\\
Department of Mathematics, George Washington University\\
and University of Maryland, College Park.\\
e-mail: przytyck@gwu.edu\\
\ \\
Witold Rosicki\\
Institute of Mathematics, Gda\'nsk University\\
e-mail: wrosicki@math.univ.gda.pl
\end{document}